\documentclass[11pt,british,letterpaper]{amsart}

\newif\ifpdf
\ifx\pdfoutput\undefined\pdffalse\else\pdfoutput=1\pdftrue\fi\newif\pdf
\ifpdf\relax\else
\usepackage[T1]{fontenc}
\fi

\usepackage{babel,eucal,url,amssymb,stmaryrd,
enumerate,amscd,
}

\usepackage[pagebackref]{hyperref}

\theoremstyle{plain}
\newtheorem{thm}{Theorem}[section]
\newtheorem{lem}[thm]{Lemma}
\newtheorem{pro}[thm]{Proposition}
\newtheorem{co}[thm]{Corollary}

\theoremstyle{definition}

\theoremstyle{remark}

\newcommand{\Gtwo}{\ifmmode{{\rm G}_2}\else{${\rm G}_2$}\fi}

\newcommand{\LC}{{\nabla^g}}

\def\sideremark#1{\ifvmode\leavevmode\fi\vadjust{\vbox to0pt{\vss
 \hbox to 0pt{\hskip\hsize\hskip1em
 \vbox{\hsize2.5cm\tiny\raggedright\pretolerance10000
 \noindent #1\hfill}\hss}\vbox to8pt{\vfil}\vss}}}%

\date{\today}
\begin{document}
\title[Twistor and Reflector Spaces of Almost Para-Quaternionic
Manifolds]%
{Twistor and Reflector Spaces of Almost Para-Quaternionic
Manifolds}
\date{\today}

\author{Stefan Ivanov}
\address[Ivanov, Minchev, Zamkovoy]{University of Sofia "St. Kl. Ohridski"\\
Faculty of Mathematics and Informatics,\\ Blvd. James Bourchier
5,\\ 1164 Sofia, Bulgaria} \email{ivanovsp@fmi.uni-sofia.bg,
minchev@fmi.uni-sofia.bg, zamkovoy@fmi.uni-sofia.bg}

\author{Ivan Minchev}

\author{Simeon Zamkovoy}

\thanks{Partially supported by a Contract 154/2005 with the University of Sofia "St. Kl.
Ohridski"}

\begin{abstract}
We investigate the integrability of natural  almost complex structures on
the twistor space of an almost para-quaternionic manifold as well
as the integrability of natural almost paracomplex structures on
the reflector space of an almost para-quaternionic manifold
constructed with the help of a para-quaternionic connection. We
show that if there is an integrable structure it is independent on
the para-quaternionic connection. In dimension four, we express
the anti-self-duality condition in terms of the Riemannian Ricci
forms with respect to the associated para-quaternionic structure.

\vspace{3mm}

Key words: almost para-quaternionic mnifolds,  anti-self-dual neutral metric,
twistor space, almost complex structures.

\vspace{3mm}

MSC: 53C15, 5350, 53C25, 53C26, 53B30
\end{abstract}

\maketitle

\setcounter{tocdepth}{2} \tableofcontents

\section{Introduction and statement of the results}

We study the geometry of structures on a differentiable manifold
related to the algebra of paraquaternions. This structure leads to
the notion of (almost) hyper-paracomplex and almost
paraquaternionic manifolds in dimensions divisible by four. These
structures are also attractive in theoretical physic since they
play a role in string theory \cite{OV,hul,Bar,Hull,C4,C5} and
integrable systems \cite{D1}. For example, hyper-paracomplex
geometry arises in connection with different versions of the c-map
\cite{C5}. New versions of the c-map are constructed in \cite{C5}
which allow the authors to obtain the target manifolds of
hypermultiplets in Euclidean theories with rigid N =2
supersymmetry. The authors show that the resulting hypermultiplet
target spaces are para-hyper-K\"ahler manifolds.

Both quaternions $H$ and paraquaternions $\tilde H$ are real
Clifford algebras, $H=C(2,0)$,  $\tilde H=C(1,1)\cong C(0,2)$. In
other words, the algebra $\tilde H$ of paraquaternions is
generated by the unity $1$ and the generators $J^0_1,J^0_2,J^0_3$
satisfying the \emph{paraquaternionic identities},
\begin{equation}\label{par1}
(J^0_1)^2=(J^0_2)^2=-(J^0_3)^2=1,\qquad
J^0_1J^0_2=-J^0_2J^0_1=J^0_3.
\end{equation}

We recall the notion of almost hyper-paracomplex manifold
introduced by Libermann \cite{Lib}. An \emph{almost quaternionic
structure of the second kind}  on a smooth manifold consists of
two almost product structures $J_1,J_2$ and an almost complex
structure $J_3$ which mutually anti-commute, i.e. these structures
satisfy the paraquaternionic identities \eqref{par1}. Such a
structure is also called \emph{complex product structure}
\cite{ASal,An}.

An \emph{almost hyper-paracomplex structure} on a 4n-dimensional
manifold $M$ is a triple $\tilde H=(J_a), a=1,2,3$, where
$J_{\alpha}$,$\alpha = 1,2$ are almost paracomplex structures
$J_{a}:TM\rightarrow TM$, and $J_{3}:TM\rightarrow TM$ is an
almost complex structure, satisfying the paraquaternionic
identities \eqref{par1}. We note that on an almost
hyper-paracomplex manifold there is actually a 1-sheeted
hyperboloid worth of almost complex structures:
$$S^2_1(-1)=\{c_1J_1+c_2J_2+c_3J_3 : c_1^2+c_2^2-c_3^2=-1\}$$ and
a 2-sheeted hyperboloid worth of almost paracomplex structures:
$$S^2_1(1)=\{b_1J_1+b_2J_2+b_3J_3 : b_1^2+b_2^2-b_3^2=1\}.$$
When
each $J_a,a=1,2,3$ is an integrable structure, $\tilde H$ is said
to be a \emph{hyper-paracomplex structure} on $M$. Such a
structure is also called sometimes \emph{pseudo-hyper-complex}
\cite{D1}.

It is well known that the structure $J_a$ is integrable if and
only if the corresponding Nijenhuis tensor $
N_a=[J_a,J_a]+J_a^2[,]- J_a[J_a,] - J_a[,J_a]$ vanishes, $N_a=0$.
In fact an almost hyper-paracomplex structure is hyper-paracomplex
if and only if any two of the three structures $J_a,a=1,2,3$ are
integrable due to the existence of a linear identity between the
three  Nijenhuis tensors \cite{IZ,BVuk0}. In this case all almost
complex structures of the two-sheeted hyperboloid  $S^2_1(-1)$ as
well as all paracomplex structures of the one-sheeted hyperboloid
$S^2_1(1)$ are integrable. Examples of hyper-paracomplex
structures on the simple Lie groups $SL(2n+1,\mathbb
R)$,$SU(n,n+1)$ are constructed in \cite{IT}.

A \emph{hyperparahermitian metric} is a pseudo Riemannian metric
which is compatible with the (almost) hyperparacomplex structure
$\tilde H=(J_a), a=1,2,3$ in the sense that the metric is
skew-symmetric with respect  to each $J_a, a=1,2,3$. Such a metric
is necessarily of neutral signature (2n,2n). Such a structure is
called \emph{(almost) hyper-paraHermitian structure}.

An \emph{almost para-quaternionic structure} on $M$ is a rank-3
subbundle ${\mathcal P} \subset End(TM)$ which is locally spanned
by an almost hyper-para-complex structure $\tilde H=(J_{a})$;
such a locally defined triple $\tilde H$ will be called admissible
basis of ${\mathcal P}$. A linear connection $\nabla$ on $TM$ is
called \emph{para-quaternionic connection} if $\nabla$ preserves
${\mathcal P}$. We denote the space all para-quaternionic
connections on an almost para-quaternionic manifold  by
$\Delta(\mathcal P)$.

An almost para-quaternionic structure is said to be a
\emph{para-quaternionic} if there is a torsion-free
para-quaternionic connection.

An almost para-quaternionic (resp. para-quaternionic) manifold
with hyperparahermitian metric is called an \emph{almost
para-quaternionic Hermitian} (resp. \emph{para-quaternionic
Hermitian}) manifold. If the Levi-Civita connection of a
para-quaternionic Hermitian manifold is para-quaternionic
connection, then the manifold is said to be
\emph{para-quaternionic K\"ahler} manifold. This condition is
equivalent to the statement that the holonomy group of $g$ is
contained in $Sp(n,{\bf R})Sp(1,{\bf R})$ for $n\ge 2$
\cite{GRio,Vuk}. A typical example is the para-quaternionic
projective space endowed with the standard para-quaternionic
K\"ahler structure \cite{BL1}. Any para-quaternionic K\"ahler
manifold of dimension $4n\ge 8$ is known to be Einstein with
scalar curvature $s$ \cite{GRio,Vuk}. If on a para-quaternionic
K\"ahler manifold there exists an admissible basis $(\tilde H)$
such that each $J_a, a=1,2,3$ is parallel with respect to the
Levi-Civita connection, then the manifold is said to be
\emph{hyper-paraK\"ahler}. Such manifolds are also called
\emph{hypersymplectic} \cite{Hit}, \emph{neutral hyper-K\"ahler}
\cite{Kam,FPed}. The equivalent characterization is that the
holonomy group of $g$ is contained in $Sp(n,{\bf R})$ if $ n\ge 2$
\cite{Vuk}.

For $n=1$ an  almost para-quaternionic structure is the same as
oriented neutral conformal structure \cite{D1,GRio,Vuk,BVuk0} and
turns out to be always quaternionic. The existence of a (local)
hyper-paracomplex structure is a strong condition since the
integrability of the (local) almost hyper-paracomplex structure
implies that the corresponding neutral conformal structure is
anti-self-dual \cite{AG,hul,IZ}.

When $n\ge 2$, the para-quaternionic condition, i.e. the existence
of torsion-free para-quaternionic connection is a strong condition
which is equivalent to the 1-integrability of the associated
$GL(n,\tilde H)$ $Sp(1,{\bf R})$ $\cong$ $GL(2n, {\bf R})$
$Sp(1,{\bf R})$- structure \cite{An,ASal}. The paraquaternionic
condition controls the Nijenhuis tensors in the sense that
$NJ_a:=N_a$ preserves the subbundle $\mathbb P$. An invariant
first order differential operator $D$ is defined on any almost
paraquaternionic manifolds which is two-step nilpotent i.e.
$D^2=0$ exactly when the structure is paraquaternionic \cite{I4}.
Paraquaternionic  structure is a type of a para-conformal
structure \cite{BE} as well as a type of generalized hypercomplex
structure \cite{Biel}.

Let $(M,\mathcal P)$ be an almost para-quaternionic manifold. The
vector bundle $\mathcal P$ carries a natural Lorentz structure of
signature (+,+,-) such that $(J_1,J_2,J_3)$ forms an orthonormal
local basis of $P$. There are two kinds of "unit sphere" bundles
according to the existence of the 1-sheeted hyperboloid $S^2_1(1)$
and the 2-sheeted hyperboloid $S^2_1(-1)$. The \emph{twistor
space} $Z^-(M)$ is the unit pseudo-sphere bundle with fibre
$S^2_1(1)$. The \emph{reflector space} $Z^+(M)$ is the unit
pseudo-sphere bundle with fibre  $S^2_1(1)$. In other words, the
fibre of $Z^-(M)$ consists of all almost complex structures
compatible with the given paraquaternionic  structure while the
fibre of $Z^+(M)$ consists of  all almost paracomplex structures
compatible with the given paraquaternionic  structure.

Keeping in mind the formal similarity with the quaternionic
geometry where there are two natural almost complex structures on
the corresponding twistor space \cite{2,11}, one observes the
existence of two naturally arising  almost complex structures
$I_1^{\nabla}, I_2^{\nabla}$ on $Z^-(M)$ and two almost
paracomplex structures $P_1^{\nabla}, P_2^{\nabla}$ on $Z^+(M)$
defined with the help of the horizontal spaces of an arbitrary
para-quaternionic connection $\nabla\in\Delta(\mathcal P)$.

The almost paracomplex structures  on the reflector space of a
4-dimensional manifold with neutral signature metric are defined
using the horizontal spaces of the Levi-Civita connection in
\cite{JR}. The authors show that one of the almost paracomplex
structure  is never integrable while the other almost paracomplex
structure is integrable if and only if the neutral metric is
anti-self-dual. The almost complex structures  on the twistor
space of a para-quaternionic K\"ahler manifold are defined and
investigated in \cite{BDM} using  the horizontal spaces of the
Levi-Civita connection. The authors show that one of the almost
complex structure is never integrable while the other almost
complex structure  is always integrable. Both construction are
generalized in the case of twistor and reflector space of a
para-quaternionic manifold in \cite{IZ}.

In the present note we investigate the dependence on the
para-quaternionic connection of these structures on the twistor
and reflector spaces over an almost para-quaternionic manifold. We
obtain conditions on the paraquaternionic connection which imply
the coincidence of the corresponding structures
(Corollary~\ref{t2.7}, Corollary~\ref{t2.72}). We show that the
existence of an integrable almost complex structure on the twistor
space (resp. the existence of an integrable almost para-complex
structure on the reflector space) does not depend on the
para-quaternionic connection and it is equivalent to the condition
that the almost para-quaternionic manifold is quaternionic
provided the dimension is bigger than four (Theorem~\ref{t2.1},
Theorem~\ref{main}).

In dimension four we find new relations between the Riemannian
Ricci forms,i.e. the 2-forms which determine the $Sp(1,\bf R
)$-component of the Riemannian curvature, which are equivalent to
the anti-self-duality of the oriented neutral conformal structure
corresponding to a given para-quaternionic structure
(Theorem~\ref{four}).

In the last section we apply our considerations to the
paraquaternionic K\"ahler manifold with torsion recently described
by the third author in  \cite{Z}.

\section{Preliminaries}

Let ${\bf \tilde H}$ be the para-quaternions and identify ${\bf
\tilde H}^{n}={\bf R}^{4n}$. To fix notation we assume that ${\bf
\tilde H}$ acts on ${\bf \tilde H}^{n}$ by right multiplication.
This defines an antihomomorphism
\begin{gather*}\lambda :\{{\rm
unit \:para-quaternions}\} =\\= \{x+j_1y+j_2z+j_3w\ |\
x^2-y^2-z^2+w^2=1\} \longrightarrow SO(2n,2n)\subset GL(4n,{\bf
R}),\end{gather*} where our convention is that $SO(2n,2n)$ acts on
${\bf \tilde H}^{n}$ on the left. Denote the image by $Sp(1,{\bf
R})$ and let $J^0_{1}= -\lambda (j_1), J^0_{2}=-\lambda (j_2),
J^0_{3}=-\lambda (j_3)$. The Lie algebra of $Sp(1,{\bf R})$ is
$sp(1,{\bf R})=span\{J^0_{1},J^0_{2},J^0_{3}\}$ and we have
\begin{equation}\nonumber
{J^0_1}^2={J^0_2}^2=-{J^0_3}^2=1,\qquad
J^0_1J^0_2=-J^0_2J^0_1=J^0_3.
\end{equation}
Define $GL(n,\tilde H)=\{A\in GL(4n,{\bf R}): A(sp(1,{\bf
R}))A^{-1}=sp(1,{\bf R}) \}$. The Lie algebra of $GL(n,\tilde H)$
is $gl(n,\tilde H)=\{A\in gl(4n,{\bf R}): AB=BA$ for all $B\in
sp(1,{\bf R})\}$.

Let $(M,\mathcal P)$ be an almost paraquaternionic manifold and
$\tilde H=(J_a), a=1,2,3$ be an admissible local basis. We shall
use the notation  $\epsilon_1=\epsilon_2=-\epsilon_3=1$.

Let $B\in \Lambda^2(TM)$. We say that $B$ is of type $(0,2)_{J_a}$
with respect to $J_a$ if $$B(J_aX,Y)=-J_aB(X,Y)$$ and denote this
space by $\Lambda_{J_a}^{0,2}$. The projection $B^{0,2}_{J_a}$ is
given by
\begin{gather*}
B^{0,2}_{J_a}(X,Y)=\frac{1}{4}\left(\epsilon_aB(X,Y)+B(J_aX,J_aY)-J_aB(J_aX,Y)-J_aB(X,J_aY)\right).
\end{gather*}
For example, the Nijenhuis tensor $N_a\in \Lambda^{0,2}_{J_a}$.

Let $\nabla\in\Delta(\mathcal P)$ be a para-quaternionic
connection on an almost paraquaternionic manifold $(M,\mathcal
P)$. This means that there exist locally defined 1-forms
$\omega_{a}, a =1,2,3$ such that
\sideremark{zzv}
\begin{gather}
 \nabla J_1=-\omega_3\otimes J_2+\omega_2\otimes J_3,\nonumber
\\\label{zzv} \nabla J_2=\omega_3\otimes J_1+\omega_1\otimes J_3,\\\nonumber \nabla
J_3=\omega_2\otimes J_1+\omega_1\otimes J_2.
\end{gather}
An easy consequence of \eqref{zzv} is that the curvature
$R^{\nabla}$ of any para-quaternionic connection
$\nabla\in\Delta(\mathcal P)$ satisfies the relations
\begin{gather}
[R^{\nabla},J_1]=-A_3\otimes J_2+A_2\otimes J_3, \qquad
A_1=d\omega_1+\omega_2\wedge\omega_3\nonumber\\
[R^{\nabla},J_2]=A_3\otimes J_1+A_1\otimes J_3, \qquad
A_2=d\omega_2+\omega_3\wedge\omega_1,\label{rel1}\\ \nonumber
[R^{\nabla},J_3]=A_2\otimes J_1+A\otimes J_2 ,
 \qquad
A_3=d\omega_3-\omega_1\wedge\omega_2.\nonumber
\end{gather}
The Ricci 2-forms of a para-quaternionic connection are defined by
\begin{gather*}
 \rho^{\nabla}_{\alpha}(X,Y)=\frac12Tr(Z\longrightarrow
J_aR^{\nabla}(X,Y)Z), \quad \alpha=1,2, \\\nonumber
\rho^{\nabla}_3(X,Y)=-\frac12Tr(Z\longrightarrow
J_3R^{\nabla}(X,Y)Z).
\end{gather*}
It is easy to see using \eqref{rel1} that the Ricci forms are
given by
$$\rho^{\nabla}_1=d\omega_1+\omega_2\wedge\omega_3, \quad
\rho^{\nabla}_2=-d\omega_2-\omega_3\wedge\omega_1, \quad
\rho^{\nabla}_3=d\omega_3-\omega_1\wedge\omega_2.$$
We split the
curvature of $\nabla$ into $gl(n,\tilde H)$-valued part
$(R^{\nabla})'$ and $sp(1,{\bf R})$-valued part $(R^{\nabla})''$
following the classical scheme (see e.g. \cite{AM,Ish,Bes})
\begin{pro}\label{p1}
The curvature of  an almost para-quaternionic
connection on $M$ splits as follows
\begin{eqnarray}\nonumber
R^{\nabla}(X,Y)&=&(R^{\nabla})'(X,Y)+{1\over
  2n}(\rho^{\nabla}_1(X,Y)J_1+\rho^{\nabla}_2(X,Y)J_2+\rho^{\nabla}_3(X,Y)J_3),\\ \nonumber & &
  [(R^{\nabla})'(X,Y),J_a]=0,\ \ \ a=1,2,3,
\end{eqnarray}
\end{pro}
Let $\Omega, \Theta$ be the curvature 2-form and the torsion
2-form of $\nabla$ on $P(M)$, respectively (see e.g. \cite{15}).
We denote the   splitting   of   the $gl(n,\tilde H)\oplus
sp(1,{\bf R})$-valued  curvature  2-form  $\Omega$  on $P(M)$
according to Proposition~\ref{p1}, by $\Omega =\Omega' +\Omega''$,
where $\Omega' $ is a $gl(n,\tilde H)$-valued
 2-form and  $\Omega''$  is  a
$sp(1,{\bf R})$-valued form.  Explicitly,  $$ \Omega ''= \Omega
''_{1}J^0_1 + \Omega''_{2}J^0_2 + \Omega''_{3}J^0_3,$$ where
$\Omega''_a, a=1,2,3$, are 2-forms.  If $\xi, \eta, \zeta \in {\bf
R^{4n}}$, then the 2-forms $\Omega''_a, a=1,2,3$, are given by
\begin{equation*}
\Omega''_a(B(\xi ),B(\eta
))=\frac{1}{2n}\rho_a(X,Y), \quad X=u(\xi), Y=u(\eta).
\end{equation*}

\section{Twistor and reflector spaces of almost para-quaternionic manifolds}

Consider the space $\tilde H_1$ of imaginary para-quaternions. It
is isomorphic to the Lorentz space $\bf R^2_1$ with a Lorentz
metric of signature (+,+,-) defined by
$<q,q'>=-Re(q\overline{q'})$, where $\overline q=-q$ is the
conjugate imaginary para-quaternion. In $\bf R^2_1$ there are two
kinds of 'unit spheres', namely the pseudo-sphere $S^2_1(1)$ of
radius $1$ (the 1-sheeted hyperboloid) which consists of all
imaginary para-quaternions of norm $1$ and the pseudo-sphere
$S^2_1(-1)$ of radius (-1) (the 2-sheeted hyperboloid) which
contains all imaginary para-quaternions of norm (-1). The
1-sheeted hyperboloid $S^2_1(1)$ carries a natural para-complex
structure while the 2-sheeted hyperboloid $S^2_1(-1)$ carries a
natural complex structure, both induced by the cross-product on
$\tilde H_1\cong \bf R^2_1$ defined by $$X \times Y=\sum_{i\not=
k}x^iy^kJ_iJ_k $$ for vectors $X=x^iJ_i,\quad Y=y^kJ_k$. Namely,
for a tangent vector $X=x^iJ_i$ to the 1-sheeted hyperboloid
$S^2_1(1)$ at a point $q_+=q_+^kJ_k$ (resp. tangent vector
$Y=y_-^kJ_k$ to the 2-sheeted hyperboloid $S^2_1(-1)$ at a point
$q_-=q_-^kJ_k$) we define $PX:=q_+\times X$ (resp. $JY=q_-\times
Y$). It is easy to check that $PX$ is again tangent vector to
$S^2_1(1)\quad and \quad P^2X=X$ (resp. $JY$ is tangent vector to
$S^2_1(-1)\quad and \quad J^2Y=-Y$).

Let $M$ be a $4n$-dimensional manifold endowed with  an almost
para-quaternionic structure $\mathcal P$. Let $J_1,J_2, J_3$ be an
admissible basis of $\mathcal P$ defined in some neighborhood of a
given point $p\in M$. Any linear frame $u$ of $T_pM$ can be
considered as an isomorphism $u:{\bf R}^{4n}\longrightarrow T_pM.$
If we pick such a frame $u$ we can define a subspace of the space
of the all endomorphisms of $T_pM$ by $u(sp(1,{\bf R}))u^{-1}.$
Clearly, this subset is a para-quaternionic structure at the point
$p$ and in the general case this para-quaternionic structure is
different from ${\mathcal P}_p.$ We define $P(M)$ to be the set of
all linear frames $u$ which satisfy $u(sp(1,{\bf
R}))u^{-1}={\mathcal P}.$ It is easy to see that $P(M)$ is a
principal frame bundle of $M$ with structure group $GL(n,\tilde
H)Sp(1,{\bf R}),$ it is also called a $GL(n,\tilde H)Sp(1,{\bf
R})$-structure on $M.$

Let $\pi :P(M)\longrightarrow M$ be the natural projection. For
each $u\in P(M)$ we consider two linear isomorphisms $j^{+}(u)$
and $j^{-}(u)$ on $T_{\pi (u)}M$ defined by
$j^{+}(u)=uJ_{1}^{0}u^{-1}$ and $j^{-}(u)=uJ_{3}^{0}u^{-1}$ . It
is easy to see that $(j(u)^{+})^{2}=id$ and $(j(u)^{-})^{2}=-id.$
For each point $p\in M$ we define $Z^{+}_{p}(M)=\{j^{+}(u): u\in
P(M), \pi (u)=p\}$ and $Z^{-}_{p}(M)=\{j^{-}(u): u\in P(M), \pi
(u)=p\}.$ In other words, $Z^{-}_{p}(M)$ is the connected
component of $J_3$ of the space of  all complex structures (resp.
$Z^{+}_{p}(M)$ is the space of  all para-complex structures) in
the tangent space $T_{p}M$ which are compatible with the almost
para-quaternionic structure on $M$.

We define the twistor space $Z^-$ of $M,$ by setting
$Z^{-}=\bigcup_{p\in M}Z^{-}_{p}(M).$ Let $H_3$ be the stabilizer
of $J^{0}_3$ in the group $GL(n,\tilde H)Sp(1,{ \bf R})$.  There
is a bijective correspondence between the symmetric space
$GL(n,\tilde H)Sp(1,{ \bf R})/H_3\cong S^2_1(-1)^+= \{ (x,y,z)\in
{\bf R}^3\ |\ x^2+y^2-z^2=-1,\ z>0\}$ and $Z^{-}_{p}(M)$ for each
$p \in M$. So we can consider $Z^{-}$ as the associated fibre
bundle of $P(M)$ with standard fibre $GL(n,\tilde H)Sp(1,{ \bf
R})/H_3.$ Hence, $P(M)$ is a principal fibre bundle over $Z^{-}$
with structure group $H_3$ and projection $j^{-}.$ We consider the
symmetric spaces $GL(n,\tilde H)Sp(1,{ \bf R})/H_3.$ We have the
following Cartan decomposition $gl(n,\tilde H)\oplus sp(1,{\bf
R})=h_3\oplus m_3$ where $$h_3=\{A\in gl(n,\tilde H)\oplus sp(1,{
\bf R}): AJ^{0}_{3}=J^0_{3}A\}$$ is the Lie algebra of $H_3$ and
$m_3=\{A\in gl(n,\tilde H)\oplus sp(1,{ \bf R}):
AJ^0_{3}=-J^0_{3}A\}.$ It is clear that $m_3$ is generated by
$J^0_{1},$ $J^0_{2},$ i.e. $m_3=span\{J^0_{1},J^0_{2} \}.$  Hence,
if $A\in m_3$ then $J^0_{3}A\in m_3.$

We proceed with defining the reflector space $Z^+$ of $M$. We put
$Z^{+}=\bigcup_{p\in M}Z^{+}_{p}(M).$ Let $H_1$ be the stabilizer
of $J^{0}_1$ in the group $GL(n,\tilde H)Sp(1,{ \bf R})$.  There
is a bijective correspondence between the symmetric space
$GL(n,\tilde H)$ $Sp(1,{ \bf R}) / H_1\cong S^2_1(1)= \{
(x,y,z)\in {\bf R}^3\ |\ x^2+y^2-z^2=1\}$ and $Z^{+}_{p}(M)$ for
each $p \in M$. So we can consider $Z^{+}$ as the associated fibre
bundle of $P(M)$ with standard fibre $GL(n,\tilde H)Sp(1,{ \bf
R})/H_1.$ Hence, $P(M)$ is a principal fibre bundle over $Z^{+}$
with structure group $H_1$ and projection $j^{+}.$ We consider the
symmetric spaces $GL(n,\tilde H)Sp(1,{ \bf R})/H_1.$ We have the
following Cartan decomposition $gl(n,\tilde H)\oplus sp(1,{\bf
R})=h_1\oplus m_1$ where $$h_1=\{A\in gl(n,\tilde H)\oplus sp(1,{
\bf R}): AJ^{0}_{1}=J^0_{1}A\}$$ is the Lie algebra of $H_1$ and
$m_1=\{A\in gl(n,\tilde H)\oplus sp(1,{ \bf R}):
AJ^0_{1}=-J^0_{1}A\}.$ It is clear that $m_1$ is generated by
$J^0_{2},$ $J^0_{3},$ i.e. $m_1=span\{J^0_{2},J^0_{3} \}.$  Hence,
if $A\in m_1$ then $J^0_{1}A\in m_1.$

Let $\nabla$ be a para-quaternionic connection on $M$, i.e.
$\nabla$ is a linear connection in the principal bundle $P(M)$
according to (\cite{15}). Note that we make no assumptions on the
torsion or on the curvature of $\nabla.$ Keeping in mind the
formal similarity with the quaternionic geometry where one uses a
quaternionic connection to define two natural almost complex
structures on the corresponding twistor space \cite{2,11,S1,S3},
we use $\nabla$ to define two almost complex structures
$I^{\nabla}_1$ and $I^{\nabla}_2$ on the twistor space $Z^-$ and
two almost para-complex structures $P^{\nabla}_1$ and
$P^{\nabla}_2$ on the reflector space $Z^+.$ Apparently, the
construction of these structures depends on the choice of the
para-quaternionic connection $\nabla.$

We denote by $A^{\ast }$ (resp.  $B(\xi )$) the fundamental vector
field (resp. the standard horizontal vector field) on $P(M)$
corresponding to $A\in gl(n,\tilde H)\oplus sp(1,{\bf R})$ (resp.
$\xi \in {\bf R^{4n}}$).

Let $u\in P(M)$ and $Q_{u}$ be the horizontal subspace of the
tangent space $T_{u}P(M)$ induced by $\nabla$ (see e.g.
\cite{15}). The vertical space i.e. the vector space tangent to a
fibre is isomorhic to $${(gl(n,\tilde H)\oplus sp(1,{\bf
R}))}^*_u=(h_3)^*_u\oplus (m_3)^*_u=(h_1)^*_u\oplus (m_1)^*_u,$$
where $(h_i)^{\ast }_{u}=\{A^{\ast }_{u}: A\in h_i\}, (m_i)^{\ast
}_{u}=\{A^{\ast }_{u}: A\in m_i\},\ i=1,3$.

Hence, $T_{u}P(M)=(h_i)^{\ast }_{u}\oplus (m_i)^{\ast }_{u}\oplus
Q_{u}$.

For each $u\in P(M)$ we put $$V^-_{j^-(u)}=j^-_{\ast u}(
(m_3)^{\ast }_{u}), H^-_{j^-(u)}=j^-_{\ast u}Q_{u} \quad
V^+_{j^+(u)}=j^+_{\ast u}( (m_1)^{\ast }_{u}),
H^+_{j(u)}=j^+_{\ast u}Q_{u}.$$ Thus we obtain  vertical and
horizontal distributions $V^-$ and $H^-$ on $Z^-$(resp. $V^+$ and
$H^+$ on $Z^+$). Since $P(M)$ is a principal fibre bundle over
$Z^-$ (resp. $Z^+$) with structure group $H_3$ (resp $H_1$) we
have $Ker j^-_{\ast u}=(h_3)^{\ast }_{u}$ (resp. $Ker j^+_{\ast
u}=(h_1)^{\ast }_{u}$).

Hence $V^-_{j^-(u)}=j^-_{\ast u}(m_3)^{\ast }_{u}$ and $j^-_{\ast
u|(m_3)^{\ast }_{u}\oplus Q_{u}}:(m_3)^{\ast }_{u}\oplus
Q_{u}\longrightarrow T_{j^-(u)}Z^-$ is an isomorphism (resp.
$V^+_{j^+(u)}=j^+_{\ast u}(m_1)^{\ast }_{u}$ and $j^+_{\ast
u|(m_1)^{\ast }_{u}\oplus Q_{u}}:(m_1)^{\ast }_{u}\oplus
Q_{u}\longrightarrow T_{j^+(u)}Z^+$ is an isomorphism).

We define two almost complex structures $I^{\nabla}_{1}$ and
$I^{\nabla}_{2}$ on $Z^-$ by
\begin{eqnarray}\label{2.2}
& &I^{\nabla}_{1}(j^-_{\ast u}A^{\ast })=j^-_{\ast
u}(J^0_{3}A)^{\ast }, \qquad I^{\nabla}_{2}(j^-_{\ast u}A^{\ast
})=-j^-_{\ast u}(J^0_{3}A)^{\ast }\\ & &I^{\nabla}_i(j^-_{\ast
u}B(\xi ))=j^-_{\ast u}B(J^0_{3}\xi ), \qquad i=1,2,\nonumber
\end{eqnarray}
for $A \in m_3, \xi \in {\bf R^{4n}}.$

Similarly,  we define two almost para-complex structures
$P^{\nabla}_1$ and $P^{\nabla}_2$ on $Z^+$ by
\begin{eqnarray}\label{2.3}
& &P^{\nabla}_{1}(j^+_{\ast u}A^{\ast })=j^+_{\ast
u}(J^0_{1}A)^{\ast }, \qquad P^{\nabla}_{2}(j^+_{\ast u}A^{\ast
})=-j^+_{\ast u}(J^0_{1}A)^{\ast }\\ & &P^{\nabla}_i(j^+_{\ast
u}B(\xi ))=j^+_{\ast u}B(J^0_{1}\xi ), \qquad i=1,2,\nonumber
\end{eqnarray}
for $A \in m_1, \xi \in {\bf R^{4n}}.$

The almost paracomplex structures \eqref{2.3} on the reflector
space of a 4-dimensional manifold with neutral signature metric
are defined using the horizontal spaces of the Levi-Civita
connection $\LC$ in \cite{JR}. The authors show that the almost
paracomplex structure $P^{\LC}_2$ is never integrable while the
almost paracomplex structure $P^{\LC}_1$ is integrable if and only
if the neutral metric is anti-self-dual. The almost complex
structures \eqref{2.2} on the twistor space of a para-quaternionic
K\"ahler manifold are defined and investigated in \cite{BDM} with
the help of  the horizontal spaces of the Levi-Civita connection.
The authors show that the almost complex structure $I^{\LC}_2$ is
never integrable while the almost complex structure $I^{\LC}_1$ is
always integrable. Both construction are generalized in the case
of twistor and reflector space of a para-quaternionic manifold in
\cite{IZ}. Twistor space of para-quaternionic K\"ahler manifold is
investigated also in \cite{DJS} where the  LeBrun's inverse
twistor  construction for quaternionic K\"ahler manifolds
\cite{LB} has been adapted to the case of para-quaternionic
K\"ahler manifolds.

We finish this section with the next useful
\begin{lem}\label{inv}
Let  $J_-\in Z^-$ be an almost complex structure
or $J_+\in Z^+$ be an almost paracomplex structure and
$B\in\Lambda^2(TM)$.

If $B^{0,2}_{J_-}=0$ for all $J_-\in Z_-$ then $B^{0,2}_{J_+}=0$
for all $J_+\in Z_+$ and vice versa.
\end{lem}
\begin{proof}
Let $J_t=\sinh tJ_1+\cosh tJ_3, t\in \mathbb R$ be an almost
complex structure in $Z_-$. Using the conditions
$B^{0,2}_{J_t}=0=B^{0,2}_{J_3}$, we calculate
\begin{gather*}
\frac12(1+\cosh 2t)B^{0,2}_{J_1} +\frac12\sinh 2t[\mathcal B]=0,
\end{gather*}
where $\mathcal B$ is a tensor field depending on $B$.

The latter leads to $B^{0,2}_{J_1}=0$. Similarly,
$B^{0,2}_{J_2}=0$ and the lemma follows.
\end{proof}

\subsection{Dependence on the  para-quaternionic connection}
In this section we investigate when different  almost
para-quaternionic connections induce the same structure on the
twistor or reflector space over an almost para-quaternionic
manofold.

Let $\nabla$ and $\nabla^{'}$ be two different almost
para-quaternionic connections on an  almost para-quaternionic
manifold $(M,\mathcal P)$. Then we have $$
\nabla^{'}_X=\nabla_X+S_X,\qquad X\in \Gamma(TM),$$ where $S_X$ is
a (1,1) tensor on $M$ and $u^{-1}(S_X)u$ belongs to $gl(n, \tilde
H)\oplus sp(1,{\bf R})$ for any $u\in P(M)$. Thus we have the
splitting
\begin{equation}\label{eq0}
S_X(Y)=S^0_X(Y)+s^1(X)J_1Y+s^2(X)J_2Y+s^3(X)J_3Y,
\end{equation}
where $X,Y\in \Gamma(TM)$, $s^i$ are 1-forms and $[S^0_X,J_i]=0,$
$i=1,2,3$.
\begin{pro}\label{t2.5}
Let $\nabla$ and $\nabla^{'}$ be two different
para-quaternionic connections on an  almost para-quaternionic
manifold $(M,\mathcal P)$. The following conditions are
equivalent:
\begin{enumerate}
\item[i).] The two almost complex structures $I_1^{\nabla}$ and
$I_1^{\nabla^{'}}$  on the twistor space $Z^-$ coincide.
\item[ii).] The 1-forms $s^1,s^2,s^3$ are related as follows
\begin{eqnarray*}
s^1(J_1X)=s^2(J_2X)=s^3(J_3X),\qquad X\in \Gamma(TM).
\end{eqnarray*}
\item[iii).] The two almost para-complex structures $P_1^{\nabla}$ and
$P_1^{\nabla^{'}}$  on the reflector space $Z^+$ coincide.
\end{enumerate}
\end{pro}
\begin{proof}
We fix a point $J$ of the twistor space $Z^-$. We have
$J=a_1J_1+a_2J_2+a_3J_3$ with $a_1^2+a_2^2-a_3^2=-1$.
Let $\pi: Z^-\longrightarrow M$ be the
natural projection and $x=\pi(J)$. The connection $\nabla$ induces
a splitting of the tangent space of $Z^-$ into vertical and
horizontal components: $T_JZ^-= V_J^- \oplus H_J^-$. Let $v$ and
$h$ be the vertical and horizontal projections corresponding to
this splitting. Let $T_JZ^-= {V^{'}}_J^- \oplus {H^{'}}_J^-$ be
the splitting induced by $\nabla^{'}$  with the projections
$v^{'}$ and $h^{'}$, respectively. It is easy to observe the
following identities
\begin{eqnarray}\label{proj} &
& v+h=1\nonumber \\ & &v^{'}+h^{'}=1 \\ & &vv^{'}=v^{'}
\nonumber\\ & & v^{'}+vh^{'}=v \nonumber
\end{eqnarray}
In fact, $V_J^-={V^{'}}_J^-$ and we may regard this space as a
subspace of  ${\mathcal P}_x$. We have that
\begin{eqnarray*}
& & V_J^-=\{W\in {\mathcal P}_x\ |\
WJ+JW=0\}=\{w_1J_1+w_2J_2+w_3J_3\ |\ w_1a_1+w_2a_2-w_3a_3=0\},
\end{eqnarray*}
where $J=a_1J_1+a_2J_2+a_3J_3$. It follows that for any
$W\in V_J^-$,
$I_1^{\nabla}(W)=I_1^{\nabla^{'}}(W)=JW$. In general, for any
$W\in T_JZ^-$, we have
\begin{eqnarray*}
& & I_1^{\nabla}(W)=J(vW)+(J\pi(W))^h \\ & &\nonumber
I_1^{\nabla^{'}}(W)=J(v^{'}W)+(J\pi(W))^{h^{'}},
\end{eqnarray*}
where
$(.)^h$ (resp. $(.)^{h^{'}}$)
denotes the horizontal lift on $Z^-$ of the corresponding vector
field on $M$ with respect to $\nabla$ (resp. $\nabla^{'}$). Using
(\ref{proj}), we calculate that
\begin{gather}\label{coin}
v(I_1^{\nabla^{'}}W)=J(v^{'}W)+v(J\pi(W))^{h^{'}}=J((v-vh^{'})W)+v(J\pi(W))^{h^{'}}=\\\nonumber
v(I_1^{\nabla}W)-J(vh^{'}W)+v(J\pi(W))^{h^{'}}. \end{gather} We
investigate the equality
\begin{eqnarray}\label{eq1}
J(vh^{'}W)=v(J\pi(W))^{h^{'}},\qquad W\in T_JZ^-.
\end{eqnarray}
Take  $W=Y^{h^{'}}, Y\in \Gamma(TM)$ in (\ref{eq1}) to get
\begin{equation}\label{eq2}
J(vY^{h^{'}})=v(JY)^{h^{'}},\qquad Y\in T_xM
\end{equation}
Hence, (\ref{eq2}) is equivalent to
$I_1^{\nabla}=I_1^{\nabla^{'}}$ because of \eqref{coin}.

Let $(U,x_1,\dots,x_{4n})$ be a local coordinate system on $M$ and
let $Y=\sum Y^i\frac{\partial}{\partial x^i}$. The horizontal lift
of $Y$ with respect to $\nabla^{'}$ at the point $J\in Z^-$ is
given by
\begin{eqnarray*}
Y^{h^{'}}_J=\sum_{i=1}^{4n}(Y^i\circ\pi)\frac{\partial}{\partial
x^i}-\sum_{s=1}^{3}a_s{\nabla^{'}}_YJ_s
\end{eqnarray*}
We calculate
\begin{eqnarray}\label{eq3}
&
&v(JY)^{h^{'}}=(JY)^{h^{'}}-h(JY)^{h^{'}}=(JY)^{h^{'}}-(JY)^{h}=\\
& &\nonumber
=\sum_{s=1}^{3}a_s(-\nabla^{'}_{JY}J_s+\nabla_{JY}J_s)=-[S_{JY},J]
\end{eqnarray}
On the other hand, we have
\begin{eqnarray}\label{eq4}
J(vY^{h^{'}})=J(Y^{h^{'}}-Y^{h})=J\sum_{s=1}^{3}a_s(-\nabla^{'}_YJ_s+\nabla_YJ_s)=-J[S_Y,J]
\end{eqnarray}

Substitute (\ref{eq3}) and (\ref{eq4}) into (\ref{eq2}) to get
that $I_1^{\nabla}=I_1^{\nabla^{'}}$ is equivalent to the
condition
\begin{eqnarray}\label{eq5}
J[S_Y,J]=[S_{JY},J], \qquad Y\in \Gamma(TM), J\in Z^-.
\end{eqnarray}
Now, \eqref{eq5} and \eqref{eq0} easily lead to the equivalence of i) and ii).

Similarly, we obtain that $P_1^{\nabla}=P_1^{\nabla^{'}}$ if and
only if
\begin{eqnarray}\label{eq5'}
P[S_Y,J]=[S_{PY},J]
\end{eqnarray}
for any choice of $P\in Z^+$ and $Y\in TM$.

The equality \eqref{eq5'} together with \eqref{eq0} implies the equivalence of ii) and iii).
\end{proof}
\begin{co}\label{t2.72}
Let $\nabla$ and $\nabla^{'}$ be two different para-quaternionic
connections on an  almost para-quaternionic manifold $(M,\mathcal
P)$. The following conditions are equivalent:
\begin{enumerate}
\item[i).] The two almost complex structures $I_2^{\nabla}$ and
$I_2^{\nabla^{'}}$  on the twistor space $Z^-$ coincide.
\item[ii).] The 1-forms $s^1,s^2,s^3$ vanish, $s_1=s_2=s_3=0$.
\item[iii).] The two almost para-complex structures $P_2^{\nabla}$ and
$P_2^{\nabla^{'}}$  on the reflector space $Z^+$ coincide.
\end{enumerate}
\end{co}
\begin{proof}
It is sufficient to observe from the proof of
Proposition~\ref{t2.5} that $I_2^{\nabla}=I_2^{\nabla^{'}}$ is
equivalent to $J[S_Y,J]=-[S_{JY},J], \quad Y\in \Gamma(TM), J\in
Z^-$ while $P_1^{\nabla}=P_1^{\nabla^{'}}$ if and only if
$P[S_Y,J]=-[S_{PY},J]$ for any choice of $P\in Z^+$ and $Y\in TM$.
Each one of the last two conditions imply $s_1=s_2=s_3=0$.
\end{proof}
\begin{co}\label{t2.7}
Let $\nabla$ and $\nabla^{'}$ be two different
para-quaternionic connections with torsion tensors
$T^{\nabla^{'}}$ and $T^{\nabla}$, respectively, on an  almost
para-quaternionic manifold $(M,\mathcal P)$. The following
conditions are equivalent:
\begin{enumerate}
\item[i).] The two almost complex structures $I_1^{\nabla}$ and
$I_1^{\nabla^{'}}$  on the twistor space $Z^-$ coincide.
\item[ii).] The $(0,2)_J$ part with respect to all $J\in\mathcal P$ of the torsion $T^{\nabla}$ and
$T^{\nabla^{'}}$ coincides,
$(T^{\nabla})^{0,2}_J=(T^{\nabla^{'}})^{0,2}_J$.
\item[iii).] The two almost para-complex structures $P_1^{\nabla}$ and
$P_1^{\nabla^{'}}$  on the reflector space $Z^+$ coincide.
\end{enumerate}
\end{co}
\begin{proof} The equivalence of i) and iii) has been proved in
Proposition~\ref{t2.7}.

Let $S=\nabla^{'}-\nabla$. Then we have
\begin{equation}\label{tr1}
T^{\nabla^{'}}(X,Y)=T^{\nabla}(X,Y)+S_X(Y)-S_Y(X).
\end{equation}
The $(0,2)_J$-part with respect to $J$ of \eqref{tr1} gives
\begin{gather}\label{eq7}
(T^{\nabla^{'}})^{0,2}_J-(T^{\nabla})^{0,2}_J=[S_{JX},J]Y-J[S_X,J]Y-[S_{JY},J]X +
J[S_Y,J]X.
\end{gather}
Suppose iii) holds. Substitute \eqref{eq5'} into the right hand
side of \eqref{eq7} and use Lemma~\ref{inv} to get
$(T^{\nabla^{'}})^{0,2}_J=(T^{\nabla})^{0,2}_J$, i.e. ii) is true.

For the converse, put  $J=J_2$ in (\ref{eq7}) and use the splitting
(\ref{eq0}) to obtain
\begin{gather}\label{eq8}
\frac12\left(T^{\nabla^{'}})^{0,2}_{J_2}-(T^{\nabla})^{0,2}_{J_2}\right)=
\left[s_1(X)+s_3(J_2X)\right]J_1Y+\left[s_1(J_2X)+s_3(X)\right]J_3Y\\\nonumber
-\left[s_1(Y)+s_3(J_2Y)\right]J_1X-\left[s_1(J_2Y)+s_3(Y)\right]J_3X.
\end{gather}
Hence,  $s_1(J_1X)=s_3(J_3X)$  is equivalent to
$(T^{\nabla^{'}})^{0,2}_{J_2}=(T^{\nabla})^{0,2}_{J_2}$
Substitute  $J=J_1$ in (\ref{eq7}) and use the splitting
(\ref{eq0}) to obtain $s_2(J_2X)=s_3(J_3X)$ is equivalent to
$(T^{\nabla^{'}})^{0,2}_{J_1}=(T^{\nabla})^{0,2}_{J_1}$. Now,
Lemma~\ref{inv} together with Proposition~\ref{t2.5} completes the
proof.
\end{proof}

\subsection{Integrability}

In this section we investigate conditions on the  para-
quaternionic connection $\nabla$ which imply the integrability of
the almost complex structure $I^{\nabla}_1$ on $Z^-$ and almost
para-complex structure $P^{\nabla}_1$ on $Z^+$. We also  show that
$I^{\nabla}_2$ and $P^{\nabla}_2$ are never  integrable i.e. for
any choice of the  para-quaternionic connection  $\nabla$ each of
these two structures has non-vanishing Nijenhuis tensor.

We denote by $IN_i, PN_i, \quad i=1,2$ the Nijenhuis tensors of
$I_i$ and $P_i$, respectively and  recall that
\begin{gather*}
IN_i(U,W)=[I_iU,I_iW]-[U,W]-I_i[I_iU,W]-I_i[U,I_iW], \quad U,W\in
\Gamma(TZ^-), \\
PN_i(U,W)=[P_iU,P_iW]+[U,W]-P_i[P_iU,W]-P_i[U,P_iW], \quad U,W\in
\Gamma(TZ^+).
\end{gather*}
\begin{pro}\label{lint}
Let $\nabla$ be a para-quaternionic connection on an almost
para-quaternionic manifold $(M,\mathcal P)$ with torsion tensor
$T^{\nabla}$.  The following conditions are equivalent:
\begin{enumerate}
\item[i).] The almost complex structure $I^{\nabla}_1$ on the twistor
space $Z^-$ of $(M,\mathcal P)$ is integrable.
\item[ii).] The $(0,2)_J$-part $(T^{\nabla})^{0,2}_J$ of the torsion with respect
to all $J\in \mathcal P$ vanishes,
\begin{equation}\label{ltor}
(T^{\nabla})^{0,2}_J=0, J\in \mathcal P
\end{equation} and  the (2,0)+(0,2) parts of the Ricci 2-forms with respect
to an admissible basis $J_1,J_2,J_3$ of $\mathcal P$ coincide in
the sense that the following identities hold
\begin{equation}\label{idric}
\rho_{a}(J_{b}X,J_{b}Y) +
\epsilon_{b}\rho_{a}(X,Y)-\epsilon_{c}\rho_{c}(J_{b}X,Y)
- \epsilon_{c}\rho_{c}(X,J_{b}Y)=0,
\end{equation}
where $\{a,b,c\}$ is a cyclic permutation of
$\{1,2,3\}$ and $\epsilon_1=\epsilon_2=-\epsilon_3=1.$
\item[iii).] The almost paracomplex structure $P^{\nabla}_1$ on the
reflector space $Z^+$ of $(M,\mathcal P)$ is integrable.
\end{enumerate}
\end{pro}
\begin{proof} Let $J_1,J_2,J_3$ be an admissible basis
 of the almost para-quaternionic structure $\mathcal
P$.

Let $hor$ be the natural projection $T_uP\longrightarrow
(m_3)^*_u\oplus Q_u,$ with $ker(hor)=(h_3)^*_u$. We define a
tensor field $I^{'}_1$ on $P(M)$ by
\begin{eqnarray*}\nonumber
& &I^{'}_{1}(U) \in  (m_3)^*_u\oplus Q_u, \\ &
&(j^-)_{*u}(I^{'}_1(U))=I_1((j^-)_{*u}U), \qquad U\in
T_uP.\nonumber
\end{eqnarray*}
For any $U, W\in \Gamma(TP(M))$ we define
$$IN^{'}_1(U,W)=hor[I^{'}_1U,I^{'}_1W]-hor[hor U, hor W]-
I^{'}_1[I^{'}_1U,hor W]-I^{'}_1[hor U,I^{'}_1W]$$

It is easy to check that $IN^{'}_1$ is a  a tensor field on
$P(M)$. We also observe that
\begin{equation}\label{in1}
j^-_{*u}(IN^{'}_1(U,W))=IN_1(j^-_{*u}U,j^-_{*u}W), \qquad U,W\in
T_uP(M)
\end{equation}
Let $A,B\in m_3$ and $\xi, \eta\in {\bf R}^{4n}$. Using the well
known general commutation relations  among  the fundamental vector
fields and standard horizontal vector fields on the principal
bundle $P(M)$ (see e.g. \cite{15}), we calculate taking into
account \eqref{in1} that
\begin{eqnarray}
& &IN_1(j^-_{*u}(A^*_u),j^-_{*u}(B^*_u))=0. \nonumber\\ &
&IN_1(j^-_{*u}(A^*_u),j^-_{*u}(B(\xi)_u))= 0.\nonumber\\\label{n1}
& &[IN_1(j^-_{*u}(B(\xi)_u)),j^-_{*u}(B(\eta)_u))]_{H^-} =
\\ & &\qquad j^-_{*u}
(B(-\Theta(B(J_3^0\xi),B(J_3^0\eta))+\Theta(B(\xi),B(\eta))\nonumber\\
& &\qquad
+J_3^0\Theta(B(J_3^0\xi),B(\eta))+J_3^0\Theta(B(\xi),B(J_3^0\eta)))_u).\nonumber\\\label{n2}
& &[IN_1(j^-_{*u}(B(\xi)_u)),j^-_{*u}(B(\eta)_u))]_{V^-}=\\&
&\qquad \{-\rho_1(B(J_3^0\xi),B(J_3^0\eta)) +
\rho_1(B(\xi),B(\eta) )\nonumber\\ & &\qquad
+\rho_2(B(J_3^0\xi),B(\eta)) + \rho_2(B(\xi),B(J_3^0\eta))\}
j^-_{*u}(J_1^0)\nonumber\\ & &\nonumber\qquad
+\{-\rho_2(B(J_3^0\xi),B(J_3^0\eta)) + \rho_2(B(\xi),B(\eta))\\ &
&\qquad -\rho_1(B(J_3^0\xi),B(\eta)) -
\rho_1(B(\xi),B(J_3^0\eta)))\} j^-_{*u}(J_2^0).\nonumber\\& &
IN_2(j^-_{*u}(A^*_u),j^-_{*u}(B(\xi)_u))=
-4j^-_{*u}(B(A\xi)_u)\not=0.\label{n33}
\end{eqnarray}
Concerning the  reflector space, let $horr$ be the natural
projection $T_uP\longrightarrow (m_1)^*_u\oplus Q_u,$ with
$ker(horr)=(h_1)^*_u$. In a very similar way as above, we
calculate
\begin{eqnarray}
& &PN_1(j^-_{*u}(A^*_u),j^-_{*u}(B^*_u))=0. \nonumber\\ &
&PN_1(j^-_{*u}(A^*_u),j^-_{*u}(B(\xi)_u))= 0\nonumber\\\label{n3}
& &[PN_1(j^-_{*u}(B(\xi)_u)),j^-_{*u}(B(\eta)_u))]_{H^-}=\\ &
&\qquad j^-_{*u}
(B(-\Theta(B(J_1^0\xi),B(J_1^0\eta))-\Theta(B(\xi),B(\eta))\nonumber\\
& &\qquad
+J_1^0\Theta(B(J_1^0\xi),B(\eta))+J_1^0\Theta(B(\xi),B(J_1^0\eta)))_u).\nonumber\\\label{n4}
& &[PN_1(j^-_{*u}(B(\xi)_u)),j^-_{*u}(B(\eta)_u))]_{V^-}=\\ &
&\qquad \{-\rho_2(B(J_1^0\xi),B(J_1^0\eta)) -
\rho_2(B(\xi),B(\eta))\nonumber\\ & &\qquad
+\rho_3(B(J_1^0\xi),B(\eta)) + \rho_3(B(\xi),B(J_1^0\eta))\}
j^-_{*u}(J_2^0)\nonumber\\ & &\nonumber\qquad
+\{-\rho_3(B(J_1^0\xi),B(J_1^0\eta)) -
\rho_3(B(\xi),B(\eta))\nonumber\\ & &\qquad
+\rho_2(B(J_1^0\xi),B(\eta)) + \rho_2(B(\xi),B(J_1^0\eta))\}
j^-_{*u}(J_3^0).\nonumber\\& &\qquad
PN_2(j^+_{*u}(A^*_u),j^+_{*u}(B(\xi)_u))=
4j^+_{*u}(B(A\xi)_u)\not=0.\label{n22}
\end{eqnarray}

Take $X=u(\xi), Y=u(\eta)$ we see that \eqref{n1}, \eqref{n2},
\eqref{n3} and \eqref{n4} are equivalent to
\begin{gather}\label{tor}
(T^{\nabla})^{0,2}_{J_3}=T^{\nabla}(J_3X,J_3Y)-T^{\nabla}(X,Y)-J_3T^{\nabla}(J_3X,Y)-J_3T^{\nabla}(X,J_3Y)=0,\\
\label{curv}
\rho^{\nabla}_1(J_3X,J_3Y)-\rho^{\nabla}_1(X,Y)-\rho^{\nabla}_2(J_3X,Y)-\rho^{\nabla}_2(X,J_3Y)=0,\\
\label{tor1}
(T^{\nabla})^{0,2}_{J_1}=T^{\nabla}(J_1X,J_1Y)+T^{\nabla}(X,Y)-J_1T^{\nabla}(J_1X,Y)-J_1T^{\nabla}(X,J_1Y)=0,\\
\label{curv1}
\rho^{\nabla}_3(J_1X,J_1Y)+\rho^{\nabla}_3(X,Y)-\rho^{\nabla}_2(J_1X,Y)-\rho^{\nabla}_2(X,J_1Y)=0,
\end{gather}
respectively.

With the help of Lemma~\ref{inv}, we see that \eqref{tor} as well
as  \eqref{tor1} is equivalent to the statement
$(T^{\nabla})^{0,2}_{J}=0$ for all local $J\in\mathcal P$. To
complete the proof we observe that each of the equalities
\eqref{curv} and \eqref{curv1} is equivalent to \eqref{idric}.
\end{proof}
The equations \eqref{n33} and \eqref{n22} in the proof of
Proposition~\ref{lint} yield
\begin{co}
Let $\nabla$ be a para-quaternionic connection on an almost
para-quaternionic manifold $(M,\mathcal P)$ with torsion tensor
$T^{\nabla}$.
\begin{enumerate}
\item The almost complex structure $I^{\nabla}_2$ on the twistor
space $Z^-$ of $(M,\mathcal P)$ is never integrable.
\item The almost paracomplex structure $P^{\nabla}_2$ on the
reflector space $Z^+$ of $(M,\mathcal P)$ is never integrable.
\end{enumerate}
\end{co}
In the 4-dimensional case we derive
\begin{thm}\label{four}
Let $(M^4,g)$ be a 4-dimensional
pseudo-Riemannian manifold with neutral metric $g$ and let
$\mathcal P$ be the para-quaternionic structure corresponding to
the conformal class generated by $g$ with a local basis
$J_1,J_2,J_3$. Then the following conditions are equivalent
\begin{enumerate}
\item[i).] The neutral metric $g$ is anti-self-dual.
\item[ii).] The Ricci forms $\rho_a^g$ of the Levi-Civita
connection $\LC$ satisfy \eqref{idric}, i.e.
\begin{equation*} \rho^g_{a}(J_{b}X,J_{b}Y) +
\epsilon_{b}\rho^g_{a}(X,Y)-\epsilon_{c}\rho^g_{c}(J_{b}X,Y)
- \epsilon_{c}\rho^g_{c}(X,J_{b}Y)=0.
\end{equation*}
\item[iii).]  The torsion condition
\eqref{ltor}  for a linear connection $\nabla$ always implies the
curvature condition \eqref{idric}.
\end{enumerate}
\end{thm}
\begin{proof}
The proof is a direct consequence of Proposition~\ref{lint},
Corolarry~\ref{t2.7} and the result in \cite{JR} (resp.
\cite{BDM}) which states that the almost para-complex structure
$P_1^{\LC}$ (resp. the almost complex structure $I_1^{\LC}$) is
integrable exactly when the neutral conformal structure generetaed
by $g$ is anti-self-dual.
\end{proof}
In higher dimensions, the curvature condition \eqref{idric} is a
consequence of the torsion condition \eqref{ltor} in the sense of
the next
\begin{thm}\label{t2.1}
Let $\nabla$ be a para-quaternionic connection
on an almost para-quaternionic 4n-dimensional $n\ge 2$ manifold
$(M,\mathcal P)$ with torsion tensor $T^{\nabla}$. Then the
following conditions are equivalent:
\begin{enumerate}
\item[i).] The almost complex structure $I^{\nabla}_1$ on the twistor
space $Z^-$ of $(M,\mathcal P)$ is integrable.
\item[ii).] The $(0,2)_J$-part $(T^{\nabla})^{0,2}_J$ of the torsion with respect
to all $J\in \mathcal P$ vanishes,
\\\centerline{$(T^{\nabla})^{0,2}_J=0, J\in \mathcal P$.}
\item[iii).] The almost paracomplex structure $P^{\nabla}_1$ on the
reflector space $Z^+$ of $(M,\mathcal P)$ is integrable.
\end{enumerate}
\end{thm}
\begin{proof} We use Proposition~\ref{lint}.
Since the connection $\nabla$ is a para-quaternionic connection,
$\nabla\in\Delta(\mathcal P)$, the  condition \eqref{ltor} yields
the next expression of the Nijenhuis tensor $N_J$ of any local
$J\in\mathcal P$,
\begin{equation}\label{pq0}
N_J(X,Y)\in
span\{J_1X,J_1Y,J_2X,J_2Y,J_3X,J_3Y\},
\end{equation}
where $J_1,J_2,J_3$ is an admissible local basis of $\mathcal P$.

To prove that ii) implies the integrability of $I_1^{\nabla}$ and
$P_1^{\nabla}$, we apply the result of Zamkovoy \cite{Z} which
states that an almost para-quaternionic 4n-maniofold $(n\ge 2)$ is
para-quaternionic if and only if the three Nijenhuis tensors
$N_1,N_2,N_3$ satisfy the condition
\begin{equation}\label{pq}
(N_1(X,Y)+N_2(X,Y)-N_3(X,Y))\in
span\{J_1X,J_1Y,J_2X,J_2Y,J_3X,J_3Y\}.
\end{equation}

Clearly, \eqref{pq} follows from \eqref{pq0} which shows that the
almost para-quaternionic 4n-manifold ($n\ge2$) $(M,\mathcal P)$ is
a para-quaternionic manifold. Let  $\nabla^0$ be a torsion-free
para-quaternionic  connection on $(M,\mathcal P)$. Then the almost
complex structure $I_1^{\nabla^0}$ on the twistor space $Z^-$ as
well as the almost paracomplex structure $P_1^{\nabla^0}$ on the
reflector space $Z^+$ are integrable \cite{IZ} and
$I_1^{\nabla}=I_1^{\nabla^0}, \quad P_1^{\nabla}=P_1^{\nabla^0}$
due to Corolarry~\ref{t2.7}.

Hence, the equivalence between i), ii) and iii) is established,
which completes the proof.
\end{proof}
From the proof of Proposition~\ref{lint} and Theorem~\ref{t2.1},
we easily derive
\begin{co}\label{cur}
Let $\nabla$ be a para-quaternionic connection on an
$4n$-dimensional ($n\ge 2$) almost para-quaternionic manifold
$(M,\mathcal P)$ with torsion tensor $T^{\nabla}$. Then the
torsion condition \eqref{ltor} implies the curvature condition
\eqref{idric}.
\end{co}

We note that Corollary~\ref{cur} generalizes the same statement
proved in the case of PQKT-connection (see below) in \cite{Z}
using the first Bianchi identity.

Theorem~\ref{t2.1} and Corollary~\ref{t2.7}  imply
\begin{co}\label{t2.8}
Let $(M,\mathcal P)$ be an almost
para-quaternionic manifold.
\begin{enumerate}
\item[i).] Among the all almost complex structures $I_1^{\nabla},
\nabla\in\Delta(\mathcal P)$ on the twistor space $Z^-$ at most
one is integrable.
\item[ii).] Among the all almost para-complex structures $P_1^{\nabla},
\nabla\in\Delta(\mathcal P)$ on the reflector space $Z^+$ at most
one is integrable.
\end{enumerate}
\end{co}

The proof  of the next theorem follows directly from the proof of
Theorem~\ref{t2.1}, Theorem~\ref{four} and Corollary~\ref{t2.8}.
\begin{thm}\label{main}
Let $(M,\mathcal P)$ be an almost para-quaternionic $4n$-manifold.
The next three conditions are equivalent:
\begin{enumerate}
\item[1).]  Either $(M,\mathcal P)$ is a para-quaternionic manifold (if
$n\ge 2$) or $(M,\mathcal P=[g])$ is anti-self dual for $n=1$.
\item[2).] There exists an integrable almost complex structure
$I^{\nabla}_1$ on the twistor space $Z^-$ which does not depend on
the para-quaternionic connection $\nabla$.
\item[3).] There exists an integrable almost paracomplex structure
$P^{\nabla}_1$ on the reflector space $Z^+$ which does not depend
on the para-quaternionic connection $\nabla$.
\end{enumerate}
\end{thm}

\section{Para-quaternionic K\"ahler manifolds with torsion }
An almost para-quaternionic Hermitian manifold $(M,\mathcal P,g)$
is called para-quaternionic K\"ahler with torsion (PQKT) if there
exists a an almost para-quaternionic Hermitian connection
$\nabla^T\in\Delta({\mathcal P})$ whose torsion tensor $T$ is a
3-form which is (1,2)+(2,1) with respect to each $J_a$, i.e. the
tensor $T(X,Y,Z):=g(T(X,Y),Z)$ is totally skew-symmetric and
satisfies the conditions
\begin{gather*}
T(X,Y,Z)=-T(J_{\alpha}X,J_{\alpha}Y,Z) -
T(J_{\alpha}X,Y,J_{\alpha}Z) -T(X,J_{\alpha}Y,J_{\alpha}Z), \quad
\alpha =1,2;\\ T(X,Y,Z)=T(J_3X,J_3Y,Z) + T(J_3X,Y,J_3Z)
+T(X,J_3Y,J_3Z).
\end{gather*}

We recall that each PQKT is a quaternionic manifold \cite{Z}. The
condition on the torsion implies that the (0,2)-part of the
torsion of a PQKT connection vanishes. Applying
Theorem~\ref{t2.1}, we obtain
\begin{thm}\label{qkt}
Let $(M,\mathcal P,\nabla^T)$ be a PQKT and
$\nabla^0\in\Delta({\mathcal P})$ be a torsion-free
para-quaternionic connection. Then
\begin{enumerate}
\item[i).] The almost complex structure $I_1^{\nabla^T}$ on the
twistor space $Z^-$ is integrable and therefore it coincides with
$I_1^{\nabla^0}$.
\item[ii).] The almost paracomplex structure $P_1^{\nabla^T}$ on the
reflector space $Z^+$ is integrable and therefore it coincides
with $P_1^{\nabla^0}$.
\end{enumerate}
\end{thm}

\bibliographystyle{hamsplain}

\providecommand{\bysame}{\leavevmode\hbox
to3em{\hrulefill}\thinspace}

\end{document}